\input amstex.tex
\input amsppt.sty   
\magnification 1200
\vsize = 9.2 true in
\hsize=6.2 true in
\NoRunningHeads        
\parskip=\medskipamount
        \lineskip=2pt\baselineskip=18pt\lineskiplimit=0pt
       
        \TagsOnRight
        \NoBlackBoxes

        \topmatter
        \title
        Quasi-Periodic Solutions of  the  \\
        Schr\"odinger Equation with\\
        Arbitrary Algebraic Nonlinearities
        \endtitle
\author
         W.-M.~Wang        \endauthor        
\address
{D\'epartement de Math\'ematique, Universit\'e Paris Sud, 91405 Orsay Cedex, FRANCE}
\endaddress
        \email
{wei-min.wang\@math.u-psud.fr}
\endemail
\abstract
We present a geometric formulation of existence of time quasi-periodic solutions.
As an application, we prove the existence of  quasi-periodic solutions
of $b$ frequencies, $b\leq d+2$, in arbitrary dimension $d$ and for arbitrary non integrable algebraic nonlinearity $p$. This reflects the conservation of $d$ momenta, energy and $L^2$ norm.  In $1d$, we prove the existence of quasi-periodic solutions with arbitrary $b$ and for arbitrary $p$, solving a problem that started Hamiltonian PDE.
\endabstract

        \bigskip\bigskip
        \bigskip
        \toc
        \bigskip
        \bigskip 
        \widestnumber\head {Table of Contents}
      
        \head 1. Introduction and statement of the theorems
        \endhead
        \head 2. The first step in the Newton scheme--$\qquad\qquad$
         extraction of parameters
        \endhead
        \head 3. The second step
        \endhead
        \head 4. Proof of Theorem 1
        \endhead
        \head 5. Proof of Theorem 3 and the cubic nonlinearity
        \endhead
        \endtoc
        \endtopmatter
        \vfill\eject
        \bigskip
\document
\head{\bf 1. Introduction and statement of the theorems}\endhead
We consider the nonlinear Schr\"odinger equation on the $d$-torus $\Bbb T^d=[0, 2\pi)^d$:
$$
i\frac\partial{\partial t}u =-\Delta u+|u|^{2p}u\qquad (p\geq 1, p\in\Bbb N),\tag 1.1
$$
with periodic boundary conditions: $u(t,x)=u(t, x+2n\pi)$, $x\in [0, 2\pi)^d$ for all $n\in\Bbb Z^d$.
The integer $p$ in (1.1) is {\it arbitrary}. The monomial nonlinearity can be replaced by polynomials:
$\sum_{p\geq p_0} a_p|u|^{2p}$, $a_p\in\Bbb R$.  For simplicity of exposition, we will only do that when $d=1$
and $p_0=1$ as otherwise it is integrable.

Let $u^{(0)}$ be a solution to the linear equation: 
$$i\frac\partial{\partial t}u^{(0)} =-\Delta u^{(0)}.\tag 1.2$$
The main purpose of this paper is to formulate two {\it geometric} conditions on $u^{(0)}$, which are 
sufficient and essentially necessary for the existence of time quasi-periodic solutions close to $u^{(0)}$
to the nonlinear equation (1.1) when $u^{(0)}$ has finite Fourier support. These two conditions remain
the same when $u^{(0)}$ has large Fourier support and are likely to be relevant there as well.

We seek quasi-periodic solutions to (1.1) with $b$ frequencies in the form 
$$
u(t, x)=\sum_{(n,j)}\hat u(n, j)e^{in\cdot\omega t}e^{ij\cdot x}, \qquad (n,j)\in\Bbb Z^{b+d},\qquad \tag 1.3
$$
with $\omega\in\Bbb Z^b$ to be determined. Writing in this form, a solution $u^{(0)}$ to (1.2) 
with $b$ frequencies $\omega^{(0)}=\{j_k^2\}_{k=1}^{b}$ ($j_k\neq 0$) has Fourier support
$$\text{supp } {\hat u}^{(0)}=\{(-e_{j_k}, j_k), k=1,...,b\},\tag 1.4$$
where $e_{j_k}$ is a unit vector in $\Bbb Z^b$. Define the bi-characteristics
$$\Cal C=\{(n,j)\in\Bbb Z^{b+d}|\pm n\cdot\omega^{(0)}+j^2=0\}.\tag 1.5$$
$\Cal C$ is the solution set in the form (1.3) to (1.2) and its complex conjugate in the Fourier space.

Define 
$$\aligned &\Cal C^+=\{(n,j) |n\cdot\omega^{(0)}+j^2=0, j\neq 0\}\cup \{(n,0) |n\cdot\omega^{(0)}=0,n_1\leq 0\}, \\
&\Cal C^-=\{(n,j) |-n\cdot\omega^{(0)}+j^2=0, j\neq 0\}\cup \{(n,0) |n\cdot\omega^{(0)}=0,n_1>0\}, \\
&\Cal C^+\cap\Cal C^-=\emptyset ,\quad \Cal C^+\cup\Cal C^-=\Cal C\text{ and }\\
&\Cal C^{++}=\{(n,j)|(n,j)=(n'-n'',j'-j''), (n',j'), (n'', j'')\in\Cal C^+\},\endaligned\tag 1.6$$
the set of differences and similarly $\Cal C^{+-}$, $\Cal C^{-+}$ and $\Cal C^{--}$.

By supp, we will always mean the Fourier support, so from now on we write $\text{supp } u^{(0)}$ for 
$\text{supp } {\hat u}^{(0)}$ etc. Let 
$$\Cal S=\text{supp } u^{(0)}\cup\text{supp } {\bar u}^{(0)},\tag 1.7$$
and $$F(u^{(0)})=|u^{(0)}|^{2p}u^{(0)}$$
the error term. For any subset $S\subset\Bbb Z^{b+d}$ and $f$ a function on $\Bbb T^{b+d}$ of the form (1.3),
define the restricted sum 
$$f_S(t,x)=\sum_{(n,j)\in S} \hat f(n, j)e^{in\cdot\omega t}e^{ij\cdot x}.$$

The two conditions on $u^{(0)}$ are:

\item{(i)} non intersection
$$\{\text{supp }F(u^{(0)})\cup \text{supp }F({\bar u}^{(0)})\}\cap\{\Cal C\backslash\Cal S\}=\emptyset,$$

\item{(ii)} non spiral
$$\{ (n,0)| n\neq 0\} \notin \text{supp } \{\prod_m [|u^{(0)}|^{2p}_{\Cal C^{++}} ]^{k_m} [|u^{(0)}|^{2(p-1)}
({\bar u}^{(0)})^2]_{\Cal C^{+-}}^{s_m}
[|u^{(0)}|^{2p}_{\Cal C^{--}}]^{k'_m}[({u}^{(0)})^2|u^{(0)}|^{2(p-1)}]_{\Cal C^{-+}}^{s_m}\},$$
where $m\geq 1$, $k_m$, $k'_m\geq 0$, $s_m=0, 1$; if $s_m=0$, then either $k_m$ or $k'_m=0$.

\noindent {\it Remark.} The cubic nonlinear Schr\"odinger equations fulfill (i) in any dimension $d$ for 
{\it all} $u^{(0)}$. When $d=1$, it also fulfills (ii) for {\it all} $u^{(0)}$,  cf. sect. 5. Parameter dependent tangentially non resonant equations studied in \cite{B1, 3, EK} generically fulfill the above two conditions, as the geometry is ``random".

The first condition indicates that the error term has no additional Fourier support on the characteristic
variety $\Cal C$ aside from ${u}^{(0)}$ and ${ {\bar u}}^{(0)}$. Since to first approximation, propagation
is along the bi-characteristics, this is  a natural condition.

The second condition provides separation in the non convex $n$ direction. 
In the $j$ direction this is from convexity (cf. Lemma 2.2), fundamentally incompatible with translation invariance. This ensures that 
the lattice points on $\Cal C$ connected by the convolution generated by $u^{(0)}$ form 
disjoint sets (in $\Cal C$) of uniformly  bounded diameter.  As we will see, this is essentially necessary. 
In the special case of periodic solutions $b=1$, separation is entirely provided by convexity
as $n$ is a scalar and there is positivity. 
It is useful to note that the right side of (ii) is a subset of $\text {supp }\{Alg |u^{(0)}|^2\}$, where $Alg |u^{(0)}|^2=\{|u^{(0)}|^{2\ell}, \ell =1, 2...\}$.

The main result is 

\proclaim
{Theorem 1}
Let 
$$u^{(0)}(t, x)=\sum_{k=1}^b  a_k e^{ij_k\cdot x}e^{-i{j^2_k}t}$$
be a solution to the linear Schr\"odinger equation (1.2) and 
$a=\{a_k\}_{k=1}^b\in(0,\delta]^b=\Cal B (0,\delta)$, $j_k\neq 0$, $k=1,..., b$.
Assume $u^{(0)}$ satisfies conditions (i, ii). There exist $C$, $c>0$, such that for all
$\epsilon\in (0,1)$, there exists $\delta_0>0$ and for all $\delta\in(0,\delta_0)$ a Cantor set 
$\Cal G$ with 
$$\text{meas }\{\Cal G\cap \Cal B (0,\delta)\}/\delta^b\geq 1-C\epsilon^c.$$ 
If $a\in\Cal G$, there is a quasi-periodic solution to the nonlinear Schr\"odinger equation (1.1)
$$u(t, x)=\sum_{k=1}^b  a_k e^{ij_k\cdot x}e^{-i{\omega_k}t}+\Cal O(\delta^{2p+1}),$$
with basic frequencies $\omega=\{\omega_k\}_{k=1}^b$ satisfying 
$$\omega_k=j^2_k+\Cal O(\delta^{2p+1}),\quad k=1,...,b,\, p\geq 1.$$
The remainder $\Cal O(\delta^{2p+1})$ is in an appropriate analytic norm on $\Bbb T^{b+d}$.
\endproclaim

\noindent{\it Remark.} The above theorem also holds when there is in addition an overall phase 
$m\neq 0$, corresponding to adding $m$ to the right side of (1.1).

The following theorems are corollaries. But since they are of a general nature, we state them as 
theorems.  

\proclaim
{Theorem 2} There are quasi-periodic solutions of $b$ frequencies for the nonlinear 
Schr\"odinger equation (1.1) for $b\leq d+2$ in any dimension $d$ and for any  
nonlinearity $p$.
\endproclaim

\noindent{\it Remark.} The above theorem holds on a set of $a=\{a_k\}_{k=1}^b\in(0,\delta]^b$ 
($0<\delta\ll 1$) of positive measure. 

For $p\geq 2$ or $d\geq 3$, this is a new result on existence of quasi-periodic solutions. Basically it 
reflects the conservation of $d$-momenta, energy and $L^2$ norm. 
For quasi-periodic solutions with cubic nonlinearity $p=1$ in dimensions 
$d=1$, $2$ see \cite{KP, B1, 2, GXY}. 
(For periodic solutions $b=1$ of the wave equation, cf \cite{BB}.)

In $d\geq 2$, (1.1) is supercritical for $p$ sufficiently large, where there are no global existence results \cite{B2}. Theorem 2
constructs a type of global solutions with precise control over the Fourier coefficients.
 
\proclaim {Theorem 3} For $d=1$, there are quasi-periodic solutions of $b$ frequencies for the nonlinear 
Schr\"odinger equation (1.1) for arbitrary $b$ and arbitrary nonlinearity $p$.
\endproclaim

\noindent{\it Remark.} The above theorem holds on a set of $a=\{a_k\}_{k=1}^b\in(0,\delta]^b$ 
($0<\delta\ll 1$) of positive measure. 

For $p\geq 2$ this is a new result. The fundamental fact in $1d$ is that the zeroes of analytic functions
or polynomials are isolated points, which enables us to go beyond integrability.

\demo{Proof of Theorem 2}
Assume $b=d+2$ and $\text{supp } {u}^{(0)}=\{(-e_{j_k}, j_k), k=1,...,b\},$ $b<d+2$ follow similarly. Since the 
right side of (ii) $\subset \text {supp }\{Alg |u^{(0)}|^2\}$, it suffices to ensure $\{(n,0)|n\neq 0\}\notin
 \text {supp }\{Alg |u^{(0)}|^2\}$. 

Writing $n=\{n_k\}_{k=1}^b$ leads to 
$$ j=\sum_{k=1}^b n_kj_k=0,\qquad j_k\in\Bbb Z^d,\tag 1.8$$
which comprises of $b$ linear equations in $n_k$, $k=1,...,b$. Since the convolution is restricted to 
$\Cal C$ and the generator $|{u}^{(0)}|^2$ is real, we have the two additional equations
$$\align n\cdot\omega^{(0)}=&\sum_{k=1}^b n_kj^2_k=0\tag 1.9\\
&\sum_{k=1}^bn_k=0.\tag 1.10\endalign$$
So if the determinant of the $(d+2)\times (d+2)$ matrix:
$$\det \pmatrix 1&1&\cdots1\\ 
j_1&j_2&\cdots j_{d+2}\\
j^2_1&j^2_2&\cdots j^2_{d+2}\endpmatrix\neq 0,$$
then $n=0$ is the unique solution to (1.8-1.10) and condition (ii) is satisfied. 

To ensure (i), we note that 
$$\text{supp }F(u^{(0)})=\text{supp }(|u^{(0)}|^{2p})+\text{supp }u^{(0)},$$
where the addition is element wise. Let $(\Delta n, \Delta j)\in \text{supp }(|u^{(0)}|^{2p})$ and 
$(-e_m, j_m)\in \text{supp }u^{(0)}$. Unless $(\Delta n, \Delta j)=(-e_{k'}+e_m, j_{k'}-j_m)$, 
in which case $$(-e_m, j_m)+(\Delta n, \Delta j)=(-e_{k'}, j_{k'})\in\Cal S,$$ the function 
$$f=(-e_m+\Delta n)\cdot\omega^{(0)}+(j_m+\Delta j)^2=\Delta n\cdot\omega^{(0)}+2j_m\Delta j+\Delta j^2$$
is not a constant function. It suffices to assume $j_k(i)$ are fixed for all $k=1, 2...,b$ and $i\neq 1$. 
$f$ is a polynomial of degree $2$ in $\{j_k(1)\}_{k=1}^{b}$. The set $\phi$ where $f=0$ has codimension $1$. There are at most $b^{2p+1}$ such $f$. So for $\{j_k(1)\}_{k=1}^{b}\notin\{\phi\}$, (i) is verified. 

Clearly similar arguments hold
for $b<d+2$ and we have proven the theorem. \hfill $\square$
\enddemo
\bigskip
\noindent{\it A sketch of proof of Theorem 1.}

We write (1.1) in the Fourier space, it becomes
$$F(\hat u)= \text{diag }(n\cdot\omega+j^2)\hat u+(\hat u*\hat v)^{*p}* \hat u=0,\tag 1.11
$$
where $(n,j)\in\Bbb Z^{b+d}$ and $\hat v=\hat{\bar u}$ and $\omega\in\Bbb R^b$ is to be determined. 
Since from now on we work with (1.11), for simplicity we drop the hat and write $u$ for $\hat u$ and $v$ for $\hat v$.
We seek solutions close to the linear solution $u^{(0)}$ of $b$ frequencies, 
$\text{supp } {u}^{(0)}=\{(-e_{j_k}, j_k), k=1,...,b\},$ with frequencies
$\omega^{(0)}=\{j_k^2\}_{k=1}^{b}$ ($j_k\neq 0$) 
and small amplitudes $a=\{a_k\}_{k=1}^b$, $\Vert a\Vert\ll 1$.

We complete (1.11) by writing the equation for the complex conjugate. So we have 
$$
\cases
\text{diag }(n\cdot\omega+j^2)u+(u*v)^{*p}* u=0,\\
\text{diag }(-n\cdot\omega+j^2)v+(u*v)^{*p}* v=0,
\endcases\tag 1.12 
$$
Denote the left side of (1.12) by $F(u, v)$. We make a Lyapunov-Schmidt decomposition into the $P$-equations:
$$ F(u, v)|_{\Bbb Z^{b+d}\backslash\Cal S}=0,$$
and the $Q$-equations:
$$ F(u, v)|_{\Cal S}=0,$$
where $\Cal S$ is as defined in (1.7).
We seek solutions such that 
$u|_\Cal S=u^{(0)}$. 
The $P$-equations are infinite dimensional and determine $u$ in the complement of $\text{supp }u^{(0)}$; 
the $Q$-equations are $2b$ dimensional and determine the frequency $\omega=\{\omega_k\}_{k=1}^b$. 

This Lyapunov-Schmidt method was introduced by Craig and Wayne \cite{CW} to construct periodic solutions
for the $1d$ wave equation. It was inspired by the multiscale analysis of Fr\"ohlich and Spencer \cite{FS}.  
This construction was further developed by Bourgain to embrace the full generality of quasi-periodic solutions
and in arbitrary dimensions $d$ \cite{B1-3}. More recently, Eliasson and Kuksin \cite{EK} developed 
a KAM version of the theory in the Schr\"odinger context. All the above results, however, pertain to parameter dependent tangentially non resonant equations.

We use a Newton scheme to solve the $P$-equations, with $u^{(0)}$ as the initial approximation. The major
difference with \cite{CW, B1-3, EK} is that (1.12) is completely resonant and there are {\it no} parameters at this initial
stage. The frequency $\omega^{(0)}$ is an {\it integer} in $\Bbb Z^b$. So we need to proceed differently.

First recall the formal scheme: the first correction
$$\Delta \pmatrix u^{(1)}\\v^{(1)}\endpmatrix= \pmatrix u^{(1)}\\v^{(1)}\endpmatrix- \pmatrix u^{(0)}\\v^{(0)}\endpmatrix
=[F'( u^{(0)}, v^{(0)}]^{-1} F( u^{(0)}, v^{(0)}),$$
where $\pmatrix u^{(1)}\\v^{(1)}\endpmatrix$ is the next approximation and $F'( u^{(0)}, v^{(0)})$ is the linearized 
operator on $\ell^2(\Bbb Z^{b+d}\times \Bbb Z^{b+d})$
$$F'=D+A,\tag 1.13$$
where 
$$
D =\pmatrix \text {diag }(n\cdot\omega+j^2)&0\\ 0& \text {diag }(-n\cdot\omega+j^2)\endpmatrix\tag 1.14$$
and
$$
A=\pmatrix (p+1)(u*v)^{*p}& p(u*v)^{*p-1}*u*u\\ p(u*v)^{*p-1}*v*v& (p+1)(u*v)^{*p}\endpmatrix
\quad  (p\geq 1),\tag 1.15$$
with $\omega=\omega^{(0)}$, $u=u^{(0)}$ and $v=v^{(0)}$.

Condition (ii) ensures invertibility of $F'$, $\Vert {F'}^{-1}\Vert\leq\Cal O (\Vert a\Vert^{-2p})$ by using bounded 
diameters of the connected sets in $\Cal C$ and  making an initial excision in $a$. Combined with condition (i),
this gives 
$$\Vert \Delta u^{(1)}\Vert=\Vert \Delta v^{(1)}\Vert\leq \Cal O(\Vert a\Vert^{2p})$$
for small $\Vert a\Vert$. Inserting this into the $Q$-equations, which determine $\omega$, we achieve amplitude-frequency
modulation:
$$\aligned \Vert\Delta\omega^{(1)}\Vert &\asymp\Cal O(\Vert a\Vert^{2p})\\
\big| \det(\frac{\partial \omega^{(1)}}{\partial a})\big|&\asymp \Cal O(\Vert a\Vert^{2p})>0\endaligned$$
ensuring transversality and moreover Diophantine $\omega^{(1)}$ on a set of $a$ of positive measure.
The tangentially non resonant perturbation theory in \cite {B3} becomes available.

The first iteration is therefore the key step and is the core of the present construction, which is summarized in 
Theorem 4 at the end of section 3. We note here that if (ii) is violated, then there will be connected sets in
$\Cal C$ of infinite diameter due to translation invariance in $n$ and typically $0\in\sigma(F'(u^{(0)}, v^{(0)}))$.
So both (i, ii) are essentially necessary. (Recall (i) prevents error propagation.) (ii) is an algebraic condition
on projected loops ensuring separation in the non convex $n$ direction in the original $\Bbb Z^{b+d}$ 
lattice. This is a much weaker separation compared to the convex $j$ direction. It suffices as we only need
it for the first step. Afterwards it is taken care of by the modulated Diophantine frequency.

The present theory relies on analyzing connected sets on the bi-characteristics $\Cal C$. For the cubic nonlinearity, locally this reduces to the well known resonance structure for the partial Birkhoff normal form \cite{KP, B1, 2, GXY}.
So cubic nonlinearity is contained in the present theory as a special case. In section 5, we show this correspondence as well as proving Theorem 3. 

This paper is the nonlinear component of the resonant perturbation theory that we are developing. For the linear theory, see \cite{W1, 2}. 
\bigskip 
 \head {\bf 2. The first step in the Newton scheme-- extraction of parameters} \endhead
 Since we seek solutions with small amplitude $a=\{a_k\}_{k=1}^b$, it is convenient to rescale:
 $a_k\to\delta^{1/2p}a_k$ ($0<\delta\ll 1$). We therefore solve instead:
 $$
\cases
\text{diag }(n\cdot\omega+j^2)u+\delta (u*v)^{*p}* u=0,\\
\text{diag }(-n\cdot\omega+j^2)v+\delta (u*v)^{*p}* v=0,
\endcases\tag 2.1
$$
with $u|_{\text{supp } u^{(0)}}=u^{(0)}=a\in (0,1]^b=\Cal B(0,1)$ and the same for $v$. The frequency $\omega=\omega(a)$ is to be determined by the $Q$-equations with the initial approximation $\omega^{(0)}=\{j_k^2\}_{k=1}^{b}$ ($j_k\neq 0$).  

We use the usual procedure of iteratively solving the $P$ and then the $Q$-equations. The main novelty here
is the first step, where $\omega=\omega^{(0)}$ is a fixed integer (vector), the opposite of a Diophantine vector.
So the usual Diophantine separation is not available. We overcome the difficulty by using the algebraic
condition (ii) to analyze the geometry of the resonant structure. After the first step, with appropriate restrictions
on $a$, $\omega=\omega(a)$ becomes Diophantine from amplitude-frequency modulation and the 
non resonant mechanism in \cite{B3} is applicable.   

\noindent{\it The first iteration.}

Using $u^{(0)}$, $v^{(0)}$ and $\omega^{(0)}$ as the initial approximation, we first 
solve the $P$-equations. This requires an estimate on the inverse of the linearized operator 
$F'(u^{(0)}, v^{(0)})$ evaluated at $\omega=\omega^{(0)}$. This section is devoted to prove the
following.

\proclaim{Lemma 2.1} Assume $u^{(0)}$, $v^{(0)}$ satisfy conditions (i, ii). Let 
$a\in (0,1]^b=\Cal B(0,1)=\Cal B\subset\Bbb R^b\backslash\{0\}$. There exist $C$, $c>0$, such that
for all $\epsilon\in (0,1)$, there exist $\delta_0>0$ and $\Cal B_\epsilon$ with 
$$\text{meas }(\Cal B_\epsilon\cap\Cal B)<C\epsilon^c.$$ If $a\in\Cal B\backslash\Cal B_\epsilon$,
then for $\delta\in (0,\delta_0)$,
$$\Vert [F'(u^{(0)}, v^{(0)})]^{-1}\Vert\leq \Cal O(\delta^{-1})\tag 2.2$$
and there exists $\beta\in (0,1)$ such that
$$|[F'(u^{(0)}, v^{(0)})]^{-1}(x,y)|\leq \delta^{\beta|x-y|}=e^{-\beta|\log\delta||x-y|}\tag 2.3$$
for all $|x-y|>1/\beta^2$.
\endproclaim 

We need the following arithmetic fact, a proof of which is given in  \cite {B3, Lemma 19.10}, so we
do not repeat it here.

\proclaim{Lemma 2.2} Fix any large number $B$. There is a partition $\{\pi_\alpha\}$ of $\Bbb Z^d$
satisfying the properties
$$\text{diam } \pi_\alpha<B^{C_0},\quad C_0=C(d),$$
and $$|j-j'|+|j^2-{j'}^2|>B,$$ if $n\in\pi_\alpha$,  $n'\in\pi_{\alpha'}$ with $\alpha\neq\alpha'$.
\endproclaim
\smallskip
\demo{Proof of Lemma 2.1}
Let $P_\pm$ be the projection on $\Bbb Z^{b+d}$ onto $\Cal C^\pm$ defined in (1.6),
$$P=\pmatrix P_+&0\\0&P_-\endpmatrix\tag 2.4$$
on $\Bbb Z^{b+d}\times \Bbb Z^{b+d}$
and $P^c$ the projection onto the complement. The linearized operator  $F'(u^{(0)}, v^{(0)})$
is $F'=D+\delta A$ with $D$ as in (1.14) and $A$ as in (1.15). 

From the Schur complement reduction \cite{S1, 2}, $\lambda\in\sigma(F')\cap [-1/2,1/2]$ if and only if
$0\in\sigma(H)$, where 
$$H=PF'P-\lambda+PF'P^c(P^cF'P^c-\lambda)^{-1}P^cF'P\tag 2.5$$
is the effective operator acting on the bi-characteristics $\Cal C$. 
Since $\Vert P^cF'P^c\Vert>1-\Cal O(\delta)>1/2$ and $\Vert PF'P^c\Vert=\Cal O(\delta)$, 
the last term in (2.5) is of order $\delta^2$, uniformly for $\lambda\in[-1/2, 1/2]$. 

So $$H=PF'P-\lambda+\Cal O(\delta^2)$$ in $L^2$ uniformly for $\lambda\in[-1/2, 1/2]$. To obtain (2.2),
it suffices to prove $\Vert [PF'P]^{-1}\Vert\leq \Cal O(\delta^{-1})$. Since 
$PF'P=\delta PAP$, where $A$ as in (1.15), let $A'=PAP$, it is important to note here that 
$A'=A'(a)$ depends on $a$ but is {\it independent} of $\delta$. 

Let $$R>2p \text{ diam } (\text{supp } u^{(0)}).\tag 2.6$$ From Lemma 2.2, there is a partition $\{\pi_\alpha\}$ 
of $\Bbb Z^d$ such that 
$$\text{ dist }(\Cal C\cap (\Bbb Z^b\times \pi_\alpha), \Cal C\cap (\Bbb Z^b\times \pi_{\alpha'}))
>R\tag 2.7$$
for $\alpha\neq\alpha'$ and $\text{diam }\pi_\alpha<R^{C_0}$, $C_0=C_0(d)$. 

Given two sets $S_1$, $S_2\subset\Bbb Z^{b+d}$, $S_1\cap S_2=\emptyset$, we say that $S_1$
is connected to $S_2$, $S_1\sim S_2$, if there exist $s_1\in S_1$, $s_2\in S_2$ such that 
$A(s_1, s_2)\neq 0$.  From (2.6, 2.7), $$\{\Cal C\cap (\Bbb Z^b\times \pi_\alpha)\}\nsim\{\Cal C\cap (\Bbb Z^b\times \pi_{\alpha'})\},$$ if $\alpha\neq \alpha'$.
 We say that a set $S$ is connected if for all $s\in S$, $s$ is connected to $\{S\backslash s\}$, $s\sim \{S\backslash s\}$.

Assume $S\subset \{\Cal C\cap (\Bbb Z^b\times \pi_\alpha)\}$ is a connected set and 
$|S|>2R^{C_0d}$. Since $\text{diam }\pi_\alpha<R^{C_0}$, there exist $s$, $s'$, $s\neq s'$, $s, s'\in\Cal C^+$
or $s, s'\in\Cal C^-$ such that $s=(n, j)$, $s'=(n', j)$ have the same $\Bbb Z^d$ coordinate. So
$s-s'=(n-n', 0)$. Since $S$ is a connected set, this contradicts condition (ii), unless $n=n'$. 
So 
$$ |S|\leq 2R^{C_0d}.\tag 2.8$$

Combining (2.6-2.8), this gives 
$$A'(a)=\oplus \Cal A_{\alpha, \beta} (a),$$
where $\Cal A_{\alpha, \beta}$ are matrices of sizes at most $2R^{C_0d}\times 2R^{C_0d}$.
Since $A(a)$ is a convolution matrix and $A(x, y)=0$ if $|x-y|>R$, there are at most 
$\sum_{N=1}^{2R^{C_0d}}2^{RN^2}\leq 2^{R^{3C_0d}}=K$ types of $\Cal A_{\alpha,\beta}$, 
which we rename as $\Cal A_k$, $1\leq k\leq K$.

For each $k$, $\det \Cal A_k(a)=P_k(a)$ is a polynomial in $a$ of degree at most 
$2^p\cdot 2R^{C_0d}\leq\log K$. Since for all $k$, $P_k$ is a non constant function on $\Cal B=(0,1]^b$, we have that there exist $C$, $c>0$, such that for all $0<\epsilon<1$,
$$\text{ meas }\{a\in \Cal B||P_k|<\epsilon, \text{ all }k\leq K\}\leq C\epsilon^c,\tag 2.9$$
(cf. e.g., \cite{Lemma 11.4, GS} and references therein).

Since $\Vert \Cal A_k\Vert\leq \Cal O(1)$, if $|\det\Cal A_k|>\epsilon$, then 
$\Vert [\Cal A_k]^{-1}\Vert\leq\Cal O(\epsilon^{-1})$. (The exponent is $1$ because of self-adjointness.)
This proves (2.2).

To summarize, there exists a partition of $\Cal C$, $\{S_\alpha\}$,
$ |S_\alpha|\leq 2R^{C_0d}$ for all $\alpha$, such that 
$S_\alpha\nsim S_{\alpha'}$, if $\alpha\neq\alpha'$, i.e., 
for all $s\in S_\alpha$, $s'\in S_{\alpha'}$, $A(s, s')=0$. Moreover 
$$\Vert [F'_{S_\alpha}]^{-1}\Vert\leq \Cal O(\delta^{-1}), $$
where $$\aligned F'_{S_\alpha}(x, y)=&F'(x, y),\, \text{ if } x,\, y\in S_\alpha\\
=&0,\qquad\quad\text{ otherwise}.\endaligned$$
To prove (2.3), we use the resolvent expansion:
$$[F']^{-1}=[\tilde F]^{-1}+[\tilde F]^{-1}\Gamma[F']^{-1},$$
where 
$$\aligned \tilde F&=\oplus_\alpha F'_{S_\alpha}\oplus_{(n,j)\notin\Cal C}\text{ diag }(\pm n\cdot\omega^{(0)}+j^2)\\
\Gamma&=F'-\tilde F.\endaligned$$
Now $$\Vert ([\tilde F]^{-1}\Gamma)^{2m}\Vert\leq\Cal O(\delta^m)\tag 2.10$$
for $m=1, 2,...$, since $[F'_{S_\alpha}]^{-1}\Gamma[F'_{S_{\alpha'}}]^{-1}=0$,
which is obvious when $\alpha=\alpha'$ and when $\alpha\neq\alpha'$ from (ii). 
(2.2, 2.10) give (2.3) for some $0<\beta<R^{-C_0d}$. \hfill $\square$
\enddemo

Using Lemma 2.1 to solve the $P$ and then the $Q$-equations, we obtain the following 
result after the first iteration. Let 
$$\Delta u^{(1)}=u^{(1)}-u^{(0)},\, \Delta v^{(1)}=v^{(1)}-v^{(0)}, \Delta \omega^{(1)}=\omega^{(1)}-\omega^{(0)}.$$
\proclaim{Proposition 2.3}
Assume $u^{(0)}$, $v^{(0)}$ satisfy conditions (i, ii). There exist $C$, $c>0$, such that
for all $\epsilon\in (0,1)$, there exist $\delta_0>0$ and $\Cal B_\epsilon$ with 
$$\text{meas }(\Cal B_\epsilon\cap\Cal B)<C\epsilon^c.$$ If $a\in\Cal B\backslash\Cal B_\epsilon$,
then for $\delta\in (0,\delta_0)$,
$$\Vert \Delta u^{(1)}\Vert =\Vert \Delta v^{(1)}\Vert_{\ell^2(\rho)}= \Cal O(\delta),\tag 2.11$$
where $\rho$ is a weight on $\Bbb Z^{b+d}$ satisfying 
$$\aligned\rho(x)=&e^{\beta |x|},\, 0<\beta<1\, \text{ for } |x|>x_0\\ 
=&1,\qquad\qquad\quad\, \,\text{  for } |x|\leq x_0.\endaligned$$
$$\align&\Vert F(u^{(1)}, v^{(1)})\Vert_{\ell^2(\rho)}=\Cal O(\delta^{3}),\tag 2.12\\
&\Vert \Delta \omega^{(1)}\Vert\asymp O(\delta),\tag 2.13\\
&\big| \det(\frac{\partial \omega^{(1)}}{\partial a})\big|\asymp \Cal O(\delta).\tag 2.14\endalign$$
Moreover $\omega^{(1)}$ is Diophantine
$$\Vert n\cdot \omega^{(1)}\Vert_{\Bbb T}\geq\frac{\kappa\delta}{|n|^\gamma},\quad n\in\Bbb Z^b\backslash\{0\},\,\kappa>0, \gamma>2b+1,\tag 2.15$$
where $\Vert\,\Vert_{\Bbb T}$ denotes the distance to integers in $\Bbb R$, $\kappa=\kappa(\epsilon)$ and $\gamma$ are independent of $\delta$.
\endproclaim

\demo{Proof}
From the Newton scheme
$$\aligned \Delta \pmatrix u^{(1)}\\ v^{(1)}\endpmatrix&=[F'(u^{(0)}, v^{(0)})]^{-1}F(u^{(0)}, v^{(0)})\\
&=[D^{-1}+\delta {F'}^{-1}AD^{-1}]F.\endaligned$$
Since 
$$\text{supp }F(u^{(0)}, v^{(0)})\cap\{\Cal C\backslash\Cal S\}=\emptyset,$$
from (i),
$$\Vert D^{-1}F\Vert_{\ell^2}=\Cal O(\delta)$$
and 
$$\Vert {F'}^{-1}AD^{-1}F\Vert_{\ell^2}=\Cal O(1),$$
where we also used (2.2).
So $$\Vert  \Delta \pmatrix u^{(1)}\\ v^{(1)}\endpmatrix\Vert_{\ell^2}=\Cal O(\delta)$$
and $$\Vert F(u^{(1)}, v^{(1)})\Vert_{\ell^2}\leq \Cal O(\Vert F''\Vert)\Vert  \Delta \pmatrix u^{(1)}\\ v^{(1)}\endpmatrix\Vert^2_{\ell^2}=\Cal O(\delta^3).$$
Using (2.3), the above two estimates hold in weighted space as well and we obtain (2.11, 2.12). 

From the $Q$ equations, the frequency modulation 
$$\Delta \omega_j^{(1)}=\frac{1}{a_j}F(u^{(0)}, v^{(0)})(-e_j,j)$$
is a rational function in $a$ of degree $p+1$ for each $j$. So 
$\Vert  \Delta \omega^{(1)}\Vert^2$ is controlled from below by a polynomial in $a$ of degree $2(p+b)$.
Using the same argument as in (2.9), we obtain (2.13).

Similarly $\det(\frac{\partial \omega^{(1)}}{\partial a})$ is a rational function in $a$ of degree at most
$(p+1)b$, and is controlled from below by a polynomial of degree at most $pb+b+2b!-2$.
So (2.9) gives (2.14), which in turn implies the Diophantine property (2.15) for appropriate $\kappa$ as $d^b a=\frac{d^b  \omega^{(1)}}{\big|\det(\frac{\partial \omega^{(1)}}{\partial a})\big|}$ and the excised set in $\omega^{(1)}$ is of measure $\Cal O(\delta)$.
\hfill$\square$
\enddemo
\bigskip 
\head{\bf 3. The second step}\endhead
Proposition 2.3, in particular (2.15) has transformed the (tangentially) resonant perturbation into
non resonant perturbation by amplitude-frequency modulation. However due to the extra power of
$\delta$, we need to proceed differently at the initial step. Essentially we redo Lemma 2.1 at
scale $N=|\log \delta|^s$ ($s>1$) with Diophantine $\omega^{(1)}$ replacing integer $\omega^{(0)}$.
This is made possible by the fact that the resonance structure remains the same at this scale.

Define the truncated linearized operator $F'_N(u^{(1)},v^{(1)})$ evaluated at $\omega^{(1)}$ as
 $$
\cases
F'_N(u^{(1)},v^{(1)})(x,y)&=F'(u^{(1)},v^{(1)})(x,y),\quad \Vert x\Vert_\infty\leq N,\, \Vert y\Vert_\infty\leq N\\
&=0,\qquad\qquad\qquad\qquad\text{ otherwise}.
\endcases\tag 3.1
$$
We have the analogue of Lemma 2.1.
\proclaim{Lemma 3.1} Assume Proposition 2.3 holds. There exist $C$, $c>0$, such that
for all $\epsilon\in (0,1)$, there exist $\delta_0>0$ and $\Cal B_\epsilon$ with 
$$\text{meas }(\Cal B_\epsilon\cap\Cal B)<C\epsilon^c,$$
where $\Cal B=(0,1]^b$.
If $a\in\Cal B\backslash\Cal B_\epsilon$,
then for $\delta\in (0,\delta_0)$,
$$\Vert [F'_N(u^{(1)}, v^{(1)})]^{-1}\Vert\leq \Cal O(\delta^{-1-\epsilon})\tag 3.2$$
and there exists $\beta\in (0,1)$ such that
$$|[F'_N(u^{(1)}, v^{(1)})]^{-1}(x,y)|\leq \delta^{\beta|x-y|}=e^{-\beta|\log\delta||x-y|}\tag 3.3$$
for all $|x-y|>1/\beta^2$.
\endproclaim 

(3.2, 3.3) together with (2.15) put the construction in the standard non resonant form as put forward by
Bourgain \cite{B3}, cf. also \cite{BW}. (2.14) ensures transversality, in other words, continuing 
amplitude-frequency modulation. The extra power of $\delta$ in (2.15) are compensated by the 
two extra powers of $\delta$ in the remainder in (2.12). 

\demo{Proof of Lemma 3.1}
Using (2.11, 2.12), it suffices to work with $F'_N(u^{(0)}, v^{(0)})$ (evaluated at $\omega^{(1)}$). The geometric 
construction in the proof of Lemma 2.1 remains valid:
$$PF'_NP=\oplus_k\Gamma_k(a),\quad k\leq K_1(N),$$
where $\Gamma_k$ are matrices of at most sizes $2R^{C_0d}\times 2R^{C_0d}$, $R$ as in (2.6),
$$\Gamma_k=\pmatrix \text{diag }(n\cdot\Delta\omega^{(1)})&0\\0&\text{diag }(-n\cdot\Delta\omega^{(1)})\endpmatrix+
\delta\Cal A_k, $$
and $\Cal A_k$ are convolution matrices.

With the addition of the diagonal term, $\Gamma_k$ is no longer a convolution matrix. Moreover 
$|n|\leq N=|\log\delta|^s$ ($s>1$) depend on $\delta$. So we need to proceed differently because of 
uniformity considerations in estimates of type (2.9). 

Fix $$N_0=\max (100 R^{C_0d}, R^{2C_0d}).\tag 3.4$$
For a given $\Gamma_k$, define the support of $\Gamma_k$ to be 
$$\Bbb Z^{b+d}\times \Bbb Z^{b+d}\supset\text{supp }\Gamma_k=\{(x,y)|\Gamma_k(x,y)\neq 0\}.$$
For matrices $\Gamma_k$, such that $$\text{supp }\Gamma_k\cap \{[-N_0, N_0]^{b+d}\times [-N_0, N_0]^{b+d}\}\neq\emptyset,\tag 3.5$$
we proceed as in the proof of Lemma 2.1. There are at most $K_0$ (independent of $\delta$) of these matrices.
Let  $P_k=P_k(a)=\det \Gamma_k(a)$, we have that there exist $C$, $c>0$, such that for all $0<\epsilon<1$,
$$\text{ meas }\{a\in \Cal B||P_k|<\delta^{1+\epsilon}, \text{ all }k\leq K_0\}\leq C\delta^{c\epsilon}.$$
So $\Vert \Gamma_k^{-1}\Vert\leq\Cal O(\delta^{-1-\epsilon})$ for all $k\leq K_0$. It is important to note
that $C, c$ are {\it independent} of $\delta$, since $P_k/\delta$ are {\it independent} of $\delta$.

For matrices $\Gamma_k$ with $k>K_0$, 
$$\text{supp }\Gamma_k\cap \{[-N_0, N_0]^{b+d}\times [-N_0, N_0]^{b+d}\}=\emptyset$$
by definition. We use perturbation theory. For any $\Gamma_k$, fix $N'$, with $|N'|>N_0$, such that
for all 
$(n, j)\in\text{supp }\Gamma_k$, we can write
$(n, j)=(N', 0)+(n', j)$ with $|n'|\leq R^{C_0d}\leq N_0/100$.

Let $\lambda$ be an eigenvalue of $\Gamma_k$ with normalized eigenfunction $\phi$. Then
$$\aligned \lambda&=(\phi, \lambda\phi)=(\phi, [\text{diag }(n\cdot\Delta\omega^{(1)})+\delta\Cal A_k]\phi)\\
&=\sum_{(n,j)} n\cdot\Delta\omega^{(1)}|\phi_{n,j}|^2+\delta (\phi, \Cal A_k\phi)\\
&=N'\cdot \Delta\omega^{(1)}+\sum_{(n',j)} n'\cdot\Delta\omega^{(1)}|\phi_{n,j}|^2
+\delta (\phi, \Cal A_k\phi)\endaligned$$
First order eigenvalue variation gives 
$$N'\cdot\frac{\partial\lambda}{\partial\omega^{(1)}}={N'}^2+\Cal O(\frac{{N'}^2}{100}).$$
So $$\frac{N'\cdot\frac{\partial\lambda}{\partial\omega^{(1)}}}{|N'|}\geq\frac{99}{100}|N'|>1.$$
Using $\big|\det(\frac{\partial \omega^{(1)}}{\partial a})\big|\asymp\delta$
from (2.14) and taking $\lambda=0$, this gives $\Vert \Gamma_k^{-1}\Vert\leq \Cal O(\delta^{-1-\epsilon})$
for all $K_0<k\leq K_1(N)$ away from a set of $a$ of measure less than $\delta^{\epsilon/2}$
as $K_1(N)\leq\Cal O(|\log\delta|^{2s})$, $s>1$. 

The estimates on $\Gamma_k$ above and the Schur reduction (2.5) give (3.2). The pointwise estimates follow as in the proof of (2.3), since the geometry of 
the resonant structure remains the same. Here we also used (2.11).\hfill$\square$
\enddemo
\smallskip
We summarize the findings in Proposition 2.3 and Lemma 3.1 in the following amplitude-frequency 
modulation theorem. For simplicity, we use $u$ to denote both the function and its Fouier series as it should 
be clear from the context.
\proclaim
{Theorem 4}
Let 
$$u^{(0)}(t, x)=\sum_{k=1}^b  a_k e^{ij_k\cdot x}e^{-i{j^2_k}t}$$
be a solution to the linear Schr\"odinger equation (1.2) and 
$a=\{a_k\}_{k=1}^b\in(0, 1]^b=\Cal B $, $j_k\neq 0$, $k=1,..., b$ and $N=|\log\delta|^s$ ($s>1$).
If $u^{(0)}$ satisfies conditions (i, ii), then there exist $C$, $c>0$, such that for all
$\epsilon\in (0,1)$,  
there exist $\delta_0>0$ and $\Cal B_\epsilon$ with 
$$\text{meas }(\Cal B_\epsilon\cap\Cal B)<C\epsilon^c.$$ If $a\in\Cal B\backslash\Cal B_\epsilon$,
then for $\delta\in (0,\delta_0)$,
$$\Vert \Delta u^{(1)}\Vert_{\ell^2(\rho)} =\Vert \Delta v^{(1)}\Vert_{\ell^2(\rho)}= \Cal O(\delta),$$
where $\rho$ is a weight on $\Bbb Z^{b+d}$ satisfying 
$$\aligned \rho(x)=&e^{\beta |x|},\, 0<\beta<1\, \text{ for } |x|>x_0\\ 
=&1,\qquad\qquad\quad\, \,\text{  for } |x|\leq x_0.\endaligned$$
$$\aligned&\Vert F(u^{(1)}, v^{(1)})\Vert_{\ell^2(\rho)}=\Cal O(\delta^{3}),\\
&\Vert \Delta \omega^{(1)}\Vert\asymp O(\delta),\\
&\big| \det(\frac{\partial \omega^{(1)}}{\partial a})\big|\asymp \Cal O(\delta),\endaligned$$
and $\omega^{(1)}$ is Diophantine
$$\Vert n\cdot \omega^{(1)}\Vert_{\Bbb T}\geq\frac{\kappa\delta}{|n|^\gamma},\quad n\in\Bbb Z^b\backslash\{0\},\,\kappa>0, \gamma>2b+1,$$
where $\Vert\,\Vert_{\Bbb T}$ denotes the distance to integers in $\Bbb R$, $\kappa=\kappa(\epsilon)$ and 
$\gamma$ are independent of $\delta$.

Moreover $$\Vert [F'_N(u^{(1)}, v^{(1)})]^{-1}\Vert\leq \Cal O(\delta^{-1-\epsilon})\tag 3.6$$
and 
$$|[F'_N(u^{(1)}, v^{(1)})]^{-1}(x,y)|\leq \delta^{\beta|x-y|}=e^{-\beta|\log\delta||x-y|}\tag 3.7$$
for all $|x-y|>1/\beta^2$, where $F'_N$ as in (3.1) with $N=|\log\delta|^s$ ($s>1$).
\endproclaim 
\bigskip 
\head{\bf 4. Proof of Theorem 1}\endhead
Theorem 4 provides the input for the initial scale in the Newton scheme in \cite{B3}. To 
continue the iteration, we need the analogues of (3.6, 3.7) at larger scales. This is achieved as follows.

Let $T=F'$ be the linearized operator defined as in (1.13-1.15) and the restrcited operator 
$T_N=F'_N$ as in (3.1). To increase the scale from $N$ to a larger scale $N_1$, we pave the $N_1$
cubes with $N$ cubes. In the $j$ direction, this is taken care of by perturbation; while in the $n$
direction by adding an additional parameter $\theta\in\Bbb R$ and consider $T(\theta)$: 
$$
T(\theta) =\pmatrix \text {diag }(n\cdot\omega+j^2+\theta)&0\\ 0& \text {diag }(-n\cdot\omega+j^2-\theta)\endpmatrix+\delta A,$$
where $A$ as defined in (1.15). 

Let $N=|\log\delta|^s$ $(s>1)$ as in Lemma 3.1 and $T_N(\theta)=T_N(\theta; u^{(1)}, v^{(1)})$ evaluated at
$\omega^{(1)}$. We have the following estimates. 
\proclaim{Lemma 4.1} Assume $u^{(0)}$, $v^{(0)}$ satisfy conditions (i, ii) and $a\in\Cal B\backslash\Cal B_\epsilon$, the set defined in Theorem 4. 
Then
$$\Vert [T_N(\theta)]^{-1}\Vert\leq \Cal O(\delta^{-1-\epsilon})\tag 4.1$$
and there exists $\beta\in (0,1)$ such that
$$|[T_N(\theta)]^{-1}(x,y)|\leq \delta^{\beta|x-y|}=e^{-\beta|\log\delta||x-y|}\tag 4.2$$
for all $|x-y|>1/\beta^2$, away from a set $B_N(\theta)\subset\Bbb R$ with 
$$\text{meas }B_N(\theta)<\delta^{1+c\epsilon}.$$
\endproclaim 
\demo{Proof} For any $\theta\in\Bbb R$, we can always write 
$\theta=\Theta+\delta\theta'$, where $\Theta\in\Bbb Z$, $-\frac{1}{2}\leq \delta\theta'<\frac{1}{2}$.
Since $\Vert T_N\Vert\leq 2 |\log\delta|^{2s}$ ($s>1$), we may restrict $\Theta$ to $|\Theta|\leq 2|\log\delta|^{2s}+1$.

Then 
$$\aligned
T_N(\theta) =&\pmatrix \text {diag }(n\cdot\omega^{(0)}+j^2+\Theta)&0\\ 0& \text {diag }(-n\cdot\omega^{(0)}+j^2-\Theta)\endpmatrix\\
&+\delta \pmatrix \text {diag }(n\cdot\tilde\omega+\theta')&0\\ 0& \text {diag }(-n\cdot\tilde\omega-\theta')\endpmatrix
+\delta A_N,\endaligned$$
where $\omega^{(0)}\in\Bbb Z^b$, $\Theta\in\Bbb Z$, $\tilde\omega=\Delta\omega^{(1)}/\delta$ is Diophantine
and $A_N$ is the restricted $A$ as defined in (1.15, 3.1).

Let $$H=\delta\left[ \pmatrix \text {diag }(n\cdot\tilde\omega+\theta')&0\\ 0& \text {diag }(-n\cdot\tilde\omega-\theta')\endpmatrix
+A_N\right].$$
Let $P_+$ be the projection onto the set $\{(n,j) |n\cdot\omega^{(0)}+j^2+\Theta=0\}$ and 
$P_-$ the projection onto the set $\{(n,j) |-n\cdot\omega^{(0)}+j^2-\Theta=0\}$ when $\Theta\neq 0$;
when $\Theta=0$, use the definition in the proof of Lemma 2.1. Define 
$$P=\pmatrix P_+&0\\0&P_-\endpmatrix$$
and $P^c$ the projection onto the complement as before. 

We proceed using the Schur reduction as 
in the proof of Lemme 2.1 and 3.1. It suffices to estimate
$PHP=\oplus_k\Gamma_k(\theta)$
as $P^cT_NP^c$ is invertible, $\Vert (P^cT_NP^c-\lambda)^{-1}\Vert\leq 4$ uniformly in $\theta$ and
$\lambda\in[-1/4, 1/4]$.

For $\Gamma_k$ such that (3.5) is satisfied, we obtain 
$$\text{ meas }\{\theta|\Vert\Gamma^{-1}_k\Vert>\delta^{-1-\epsilon}\}\leq C\delta^{c\epsilon}\qquad (\epsilon>0),$$
via the estimates on $\det \Gamma_k$ as in the proof of Lemma 3.1.
For other $\Gamma_k$, we make a shift
$n\cdot\tilde\omega=N'\cdot \tilde\omega+n'\cdot \tilde\omega,$
where $|N'|<|\log\delta|^s$, $|n'|<N_0$, $N_0$ as in (3.4). Let $\theta''=\theta'+ N'\cdot \tilde\omega$. 
Then clearly 
$$\text{ meas }\{\theta''|\text{dist }(\theta'', \sigma(\Gamma^k))\leq \delta^{1+\epsilon}\}\leq \Cal O(\delta^{1+\epsilon}).$$

Combining the above estimates on $\Gamma_k$, we obtain (4.1). Using the geometric structure afforded by the resonance structure,
we obtain (4.2) as before.\hfill$\square$ 
\enddemo
\smallskip
\demo{Proof of Theorem 1}
This follows from Lemma 4.1 and lemme 19.13, 19.38 and 19.65 in chapter 19 of \cite{B3}. 
\hfill$\square$
\enddemo
\bigskip
\head{\bf 5. Proof of Theorem 3 and the cubic nonlinearity}\endhead
\noindent {\it Cubic nonlinearity.}

For simplicity we write $u$ for $u^{(0)}$ and $\omega$ for $\omega^{(0)}$, the solutions and frequencies
of the linear equation. The symbols of convolution for the cubic nonlinearity are $|u|^2$, $u^2$ and ${\bar u}^2$.
Assume $(n,j)\in\Cal C^+$ ($(n,j)\in\Cal C^-$ works similarly). In order that $(n,j)$ is connected 
to $(n'j')\in\Cal C$, it is necessary that either
\item{(a)} $[u*v](n,j;n'j')\neq 0$ or
\item{(b)} $[v*v](n,j;n'j')\neq 0$.

\noindent Case (a): Since $$\aligned n\cdot\omega+j^2&=0,\\
n'\cdot\omega+{j'}^2&=0,\endaligned$$
subtracting the two equations gives immediately
$$(j_k-j_k')\cdot(j+j_k)=0,\tag 5.1$$
where $j_k, j_{k'}\in\Bbb Z^d$ ($k, k'=1,..., b$), are the $b$ Fourier components of $u$.

\noindent Case (b): Since $$\aligned n\cdot\omega+j^2&=0,\\
-n'\cdot\omega+{j'}^2&=0,\endaligned$$
adding the two equations gives immediately
$$(j-j_k)\cdot(j-j_{k'})=0,\tag 5.2$$
where $j_k, j_{k'}\in\Bbb Z^d$ ($k, k'=1,..., b$), are the $b$ Fourier components of $u$.

(5.1, 5.2) are precisely the well known resonant set for the partial Birkhoff normal form
transform in \cite{KP, B1, 2, GSY}. (5.1, 5.2) describe rectangular type of geometry. 
(i) is satisfied for the cubic nonlinear Schr\"odinger in any $d$. 
In $1d$, (5.1, 5.2) reduce to a finite set of $2b$ lattice points in $\Bbb Z$: $j=\pm j_k$, $k=1,...,b$. 
The set $\{\mp e_k, \pm j_k\}$ is the only connected set in $\Cal C$ and (ii) is verified for all $u^{(0)}$. 

The following proof extends the existence of quasi-periodic solutions in $1d$ of any number of frequencies $b$ to arbitrary algebraic nonlinearities $p$. It shows that it is not necessary to know explicitly the locations of the resonances, which is the case for cubic nonlinearity, except that they are isolated. 
\smallskip
\noindent {\it Proof of Theorem 3.}

Let $$\Gamma^+=\text{supp }[|u^{(0)}|^{2p}]=\{(\Delta n, \Delta j)\},$$ 
where $$\aligned &\Delta n=-\sum p_{kk'}(e_k-e_{k'}),\\
&\Delta j=\sum p_{kk'}(j_k-j_{k'}),\\
&p_{kk'}\geq 0,\, \sum p_{kk'}\leq p,\, k, k'=1,...,b,\\
&\{(-e_k, j_k)\}_{k=1}^b=\text{supp } u^{(0)},\endaligned\tag 5.3$$
and 
$$\Gamma^-=\text{supp }[|u^{(0)}|^{2(p-1)}{v^{(0)}}^2]=\{(\Delta n, \Delta j)\},$$ 
where 
$$\aligned &\Delta n=-\sum p_{kk'}(e_k-e_{k'})+(e_{\kappa}+e_{\kappa'}),\\
&\Delta j=\sum p_{kk'}(j_k-j_{k'})-(j_{\kappa}+j_{\kappa'}),\\
&p_{kk'}\geq 0,\, \sum p_{kk'}\leq p-1,\, k, k', \kappa, \kappa'=1,...,b.\endaligned\tag 5.4$$

Assume $\{j_k\}_{k=1}^b$ is such that $(\Delta n, 0)\notin \Gamma^+\cup \Gamma^-$ for $n\neq 0$, which
is always possible from (5.3, 5.4).
Assume $(n,j)\in\Cal C^+$ ($(n,j)\in\Cal C^-$ works similarly). Then
in order for $(n,j)$ to be connected to $(n'j')\in\Cal C$ via $(\Delta n,\Delta j)$, at least one of the
following two equations has to be satisfied as in the cubic case:
$$2j\Delta j+\Delta j^2+\Delta n\cdot\omega=0,\qquad \text{if }(\Delta n,\Delta j)\in \Gamma^+\tag 5.5$$
or 
$$2j^2+2j\Delta j+\Delta j^2-\Delta n\cdot\omega=0, \qquad \text{if }(\Delta n,\Delta j)\in \Gamma^-.\tag 5.6$$

Let $\{j, j+\Delta j\}$ be the $\Bbb Z$ coordinates of the set of connected pairs determined by (5.5, 5.6).
Let $D=\{j-j',  j+\Delta j-j', j-j'-\Delta j', j+\Delta j-j'-\Delta j', \Delta j\neq \Delta j'\}$ be the 
set of differences. The constant functions in $D$ are only generated by pairs where 
$\Delta j=j_m\pm j_n$ (and similarly for $\Delta j'$), in which case $(j, j+\Delta j)=(\mp j_n,  j_m)$. These are the terms  present in the cubic nonlinearity. 

For the non constant functions $f$, $f=0$ implies $P=0$ for some polynomial $P$ in 
$\{j_k\}_{k=1}^b$. Let $\phi$ be the zero
set of $P$. It has codimension $1$. There are at most $8b^{2p}$ such polynomials. 
Choose $\{j_k\}_{k=1}^b\notin\{\phi\}$. Condition (ii) is then verified. Fulfilling (i) as in the proof of 
Theorem 2 concludes the present proof. \hfill $\square$

\noindent {\it Remark.} Theorem 3 is coherent with the result in \cite {B4} proving existence of almost periodic solutions with subexponential decay amplitudes for the parameter dependent non resonant $1d$ quintic nonlinear Schr\"odinger
equation.

\enddocument
\end